\documentstyle[11pt]{article}

\begin{document}

\begin{center}
\hrule\hrule

\bigskip

{\Large \bf Higher dimensional analogues of the map colouring problem}

\bigskip

\hrule

\medskip

{\large {\bf Bhaskar Bagchi and Basudeb Datta}\footnote{Research
partially supported by grant from UGC Centre for Advanced Study.}}

\medskip

\hrule\hrule

\smallskip
\end{center}

\begin{abstract}
{After a brief discussion of the history of the problem, we
propose a generalization of the map colouring problem to higher
dimensions.}
\end{abstract}

\bigskip

The map colouring problem was originally posed by Francis Guthrie
in 1852. It asks for the minimum number of colours required to
colour all possible maps (real or imagined) if one wishes to
ensure that neighbouring countries receive different colours.
Clearly one may imagine any number of countries sharing a common
point in their boundary. So, in this problem, two countries are to
be considered as neighbours only if their boundaries share an
entire linear continuum. Clearly, four colours are required since
there are maps with four mutually adjacent countries (as in the
usual planar drawing of a tetrahedron). Francis guessed that four
colours always suffice for his problem. This came to be known as
the four colour conjecture.

Beginning with Augustus de Morgan (to whom Francis' brother
Frederick first communicated the problem), innumerable
mathematicians - both professionals and amateurs - tried their
hands at the problem. Many possible approaches, reformulations,
partial results as well as several wrong proofs resulted. One
preliminary observation that could be made was that six colours
always suffice. To see this, one may proceed as follows. The map
to be coloured may well be taken as drawn on the surface of a ball
(as in a globe) rather than on a sheet of paper. A little
thought shows that, without loss of generality, the countries may
be taken to be the faces of a convex polyhedron inscribed in the
ball. Then we can invoke Euler's formula $v-e+f=2$, where $v$,
$e$, $f$ are the number of vertices, edges and faces of the
polyhedron. Using this formula, it is not hard to see that, in any
map, there must be a country with five or fewer neighbours. Now,
with six colours in hand, one may proceed to colour all maps
inductively, as follows. Assume we already know how to colour all
maps with fewer countries using six colours. Pick out a country
with five or fewer neighbours in the given map, and momentarily
forget it. By our assumption, the rest of the map can be coloured
using the six colours. Since the forgotten country has at most
five neighbours, at most five of the colours have been used to
colour these neighbours. So, there is a left over colour which can
now be used to colour the forgotten country, thus completing the
six - colouring of the given map. For an alternative proof of this
``six colours theorem", which avoids Euler's formula, read on!

The first significant progress on the four colour conjecture was
made by A. B. Kempe, when he published a wrong proof of the
conjecture in 1890. This may sound like sarcasm, but Kempe's work
contained a very fruitful idea which led to the eventual
resolution of the problem. In 1890 itself, P. J. Heawood pointed
out the flaw in Kempe's argument, but he also showed that Kempe's
logic may be modified to prove that five colours always suffice.
The argument begins as in the proof of the six colours theorem
sketched above. But, if the forgotten country has five neighbours
who have received all five colours in the colouring of the rest of
the map, then Kempe's argument gives a method to recolour them, so
that one of the colours is freed for use on  the forgotten
country.

Finally, in 1976, K. Appel and W. Haken collaborated with
a Cray supercomputer (taking more than a thousand hours of
computer time) to prove the four colour conjecture (see \cite{ah1,
ah2, ah3}). Their work made no use of the huge superstructure of
theories created during the twentieth century. Instead, they went
back to a vastly elaborated version of Kempe's original idea!
Perhaps there is a moral lurking here! Much has been made of the
fact that Appel and Haken's proof is not a ``human" one, and
nobody can possibly verify it manually. But we think that the real
problem lies elsewhere. If there was a prior theoretical proof
that one only needed to check a million cases (say) to settle the
question one way or the other, then, we think, most mathematicians
would have happily left these million verifications to a machine.
Unfortunately, the situation with the Appel-Haken proof is
not like that. Until the computer stopped and came out with its
verdict in favour of the four colour conjecture, there was no
guarantee   that the program would ever halt. It was rather like
the celebrated halting problem of Turing acted out in real life.
The proof consists in finding a finite ``unavoidable set of
reducible configurations". (A set of configurations is unavoidable
if there is a proof showing that any hypothetical counterexample
to the four colour conjecture must contain one of them. A
configuration is reducible if, whenever a counterexample contains
it, there is a well defined procedure to create a smaller
counterexample by eliminating it.) The program halted precisely
because it found such a set. Of course, one may ignore the
mechanical genesis of this list, and get its authenticity verified
(necessarily by another machine, since the list is so large). Such
a verification constitutes an airtight proof which ought to keep
everybody happy. But the question remains\,: why does such a set
exist\,? Is it mere happenstance, or is there a good theoretical
reason for its existence\,? So, the search for a ``human" proof is
still on.

Early in the game, mathematicians had realized that the infinite
variety of shapes and sizes of countries is not relevant to the
problem. What matters is the knowledge of which pairs of countries
are neighbours. In order to forget the inessential, they used the
notion of a graph. In simple terms, a graph is a picture
consisting of finitely many dots (called vertices) in which some
pairs of dots are joined by (possibly curved) line segments,
called edges. Obviously, an edge joining the vertices $x$ and $y$
intersects an edge joining the vertices $x$ and $z$ at least in
the common vertex $x$. A graph drawn in the plane is called planar
if the edges have no further intersections. Each country in a map
may be represented by a dot in its interior (perhaps its capital
city). Clearly, whenever two countries are neighbours, the
corresponding dots may be joined by a line segment (perhaps
curved) which meets their common boundary at a point not lying on
any other country, otherwise staying in the interior of these two
countries. A planar graph results. If one can colour the vertices
of this graph in four colours so that  neighbouring vertices
(i.e., those joined by an edge) receive different colours, then
one may transfer the colours from the vertices to the
corresponding countries, resulting in a proper four - colouring of
the map. Conversely, given any planar graph, one may suitably blow
up the vertices into countries, to create a corresponding map.
Thus, we have the mathematician's favourite version of the map
colouring problem\,: show that the vertices of any planar graph
can be properly coloured using four colours.

So what should be a 3-dimensional version of this problem? One
might try to define a ``spatial graph" as a graph which may be
drawn in three dimensional space without undue intersections.
Unfortunately, a moment's thought shows that all graphs are
``spatial". Indeed, one may take any (finite) number of points in
space ``in general position", i.e., such that no four of them are
on a common plane. Then one may join each pair of these points by
a straight line segment. No two of these segments will then meet,
except at a common vertex. The same observation goes for all
higher dimensions since the three dimensional space embeds in all
Euclidean spaces of higher dimensions.

Is that the end of the road as far as our search for a higher
dimensional analogue of the map colouring problem goes? Not quite!
After all, the attempt to abstract away the irrelevant geometric
details of shape and size by going to planar graphs is not a very
successful one. Indeed, any (sufficiently complicated) planar
graph is visually indistinguishable from a map. Consider the
following problem instead. Consider a finite set of
non-overlapping discs in the plane. That is, any two of them are
either disjoint or (externally) tangential. In other words, no two
of them have any common interior point. Question\,: what is the
minimum number of colours needed if we are to colour these discs
in such a way that touching discs are given different colours?
Does not that look like a much simpler problem which should have a
fairly straightforward answer? Notice that, given such a set of
discs, one can again form a planar graph (apparently of a very
special kind) by joining the centers of each pair of touching
discs by straight line segments. So, if Appel and Haken are
to be believed, then in this problem also the answer should be
four. (One can easily have four mutually touching discs in the
plane, say by taking three equal mutually touching discs and then
placing a fourth small disc in the nich\'{e} created by them.) An
amazing theorem due to Paul Koebe, E. M. Andreev and William P.
Thurston (see \cite{an, k, t}) says that every planar graph can be
redrawn as the graph of a set of non-overlapping discs. Thus, our
innocent looking problem of colouring discs is equivalent to the
map colouring problem! Now, it is obvious how one may reprove the
``six colours theorem" (albeit using the powerful K-A-T theorem).
Given any finite set of non-overlapping discs in the plane, choose
the smallest one. Say its radius is of unit length. Then, the
discs touching it are equal or larger, so that it is obvious that
at most six discs touch the smallest one. Moreover, if six discs
do touch the smallest disc, then these neighbouring discs must be
unit discs as well. We may assume that our set of discs is
``connected" (in the sense that a point particle may go from one
disc to any other, all the time staying within the discs in the
set). Now, is it possible that all the discs have six or more
discs touching them? If so, then the above argument shows that all
the discs must be of the same size, and each must have exactly six
neighbours touching it. But this is impossible since it is
intuitively clear that, in this case, a disc in the periphery
(technically, a disc touching the boundary of the convex hull of
all the discs) would have at most four neighbours! Thus, we have
shown that, given any finite configuration of non-overlapping
discs in the plane, there is at least one disc which touches at
most five others. Now, one may complete the proof of the ``six
colours theorem" as before. As promised, we have avoided any use of
Euler's formula.

Now, it should be clear how one can generalize the map colouring
problem to arbitrary dimensions. For $d \geq 1$, let
$\chi(d\hspace{.2mm})$ be the smallest number such that any
(finite) set of $d$-dimensional non-overlapping closed balls (not
necessarily of the same size) in $d$-dimensional Euclidean space
may be coloured in $\chi(d\hspace{.2mm})$ colours so that any two
touching balls receive different colours. Thus, we have $\chi(1) =
2$ (trivial\,!) and $\chi(2) =4$ (Appel and Haken\,!). What is the
value of $\chi(3)$\,?

There is a small but important question that needs to be answered.
How do we know that $\chi(d\hspace{.2mm})$ is finite\,? In other
words, is it possible that in some dimension there are arbitrarily
complicated sets of non-overlapping balls requiring unboundedly
large number of colours\,? No\,! To see this, recall that the
kissing number $\kappa(d\hspace{.2mm})$ of $d$-dimensional
Euclidean space is defined as the maximum number of
non-overlapping equal balls which may touch a given ball of the
same size. It is intuitively clear that $\kappa(d\hspace{.2mm})$
is finite for each $d$. Indeed, the $\kappa(d\hspace{.2mm})+1$
unit balls in such a kissing configuration are contained in a ball
of radius 3 (concentric with the central ball).
Comparing ($d$-dimensional) volumes, one then sees that
$\kappa(d\hspace{.2mm})+1 \leq 3^d$. Thus, $\kappa(d\hspace{.2mm})
\leq 3^d-1$ for all $d$. (This is a very crude bound. For improved
bounds, see Conway and Sloane \cite{cs}.) The exact value of
$\kappa(d\hspace{.2mm})$ is known only for $d=1, 2, 3, 4, 8$ and
24. Indeed, we have $\kappa(1)=2$, $\kappa(2)=6$, $\kappa(3)=12$,
$\kappa(4)=24$, $\kappa(8)=240$ and $\kappa(24)=196560$ (see
\cite{cs, l, m}). Isn't that a surprise\,?

Now, given any finite set of non-overlapping balls in
$d$-dimensional space, the neighbours of the smallest ball (say a
unit ball) may be shrank into unit balls, still
remaining non-overlapping and still touching the smallest ball.
This shows that there is at least one ball (namely the smallest)
with at most $\kappa(d\hspace{.2mm})$ neighbours. Therefore, as
before, one may prove by induction on the number of balls that any
finite set of non-overlapping balls in $d$-space may be properly coloured
using at most $\kappa(d\hspace{.2mm}) +1$ colours. Thus,
$\chi(d\hspace{.2mm}) \leq \kappa(d\hspace{.2mm})+1$. Also, in
space of dimension $d\geq 2$, one can construct a set of $d+2$
mutually touching balls (one may take $d+1$ equal balls with
centers at the vertices of a regular simplex, having the side
length of the simplex as their diameters, and then place a small
ball in the middle touching all of them). Thus,
$\chi(d\hspace{.2mm}) \geq d+2$ for all $d \geq 2$. So, we have\,:
$$
d+2 \leq \chi(d\hspace{.2mm}) \leq \kappa(d\hspace{.2mm}) +1 \leq
3^d, ~ \mbox{ for all } ~ d\geq 2.
$$

In the magic dimensions $d= 1, 2, 8$ and 24, it is known that, up
to congruence, there are unique configurations attaining the
kissing numbers. In these dimensions one may argue (exactly as we
have done in the case $d=2$) that $\chi(d\hspace{.2mm}) \leq
\kappa(d\hspace{.2mm})$, a slight improvement\,!

In our familiar three dimensional space, we have $5 \leq
\chi(3)\leq 13$. What is the true value\,? We won't even hazard a
guess. Happy hunting\,!

\bigskip

{\small
\noindent {\bf Bhaskar Bagchi} \newline
{Theoretical Statistics and Mathematics
Unit, Indian Statistical Institute,  Bangalore 560\,059, India.
E-mail\,: bbagchi@isibang.ac.in}

\smallskip

\noindent {\bf Basudeb Datta} \newline {Department of Mathematics,
Indian Institute of Science, Bangalore 560\,012,  India. \newline
E-mail\,:   dattab@math.iisc.ernet.in} }

\end{document}